\documentstyle[11pt]{article}
\begin{document}

\begin{center}
{\bf ON DEFORMATIONS AND MODULI } \\
\mbox{\ } \\
{\bf OF GENUS 2 FIBRATIONS} \\
\mbox{\ } \\
{\bf Hur\c{s}it \"{O}nsiper}
\end{center}

\begin{abstract}
In this paper, we investigate the moduli of surfaces of general
type admitting genus 2 fibrations with irregularity $q = g_b + 1$,
where
$g_b \ge 2$ is the genus of the base. We prove that
smooth fibrations are parametrized by a unique
component in the moduli space. The same result applies to nonsmooth
fibrations with special values of $g_b$. In the general case, we give
a bound
on the dimension of the corresponding connected components.
\end{abstract}

\renewcommand{\thefootnote}{}

Surfaces admitting genus 2 fibrations appear as important examples of
special surfaces of general type and there has been considerable
effort
to analyze their structure and moduli. In particular, in ([S1], [S2])
Seiler obtained substantial results on the structure of the moduli
spaces
of surfaces admitting genus 2 albanese fibrations, via detailed study
of
deformations of ruled surfaces covered by surfaces admitting genus 2
fibrations.

In this paper, we take up the same problem from a more elementary
point
of view; we apply results and methods from the theory of curves to
study
surfaces admitting genus 2 fibrations over curves of genus $g_{b} \ge
2$
with irregularity $q = g_{b} + 1$.
The basic idea is to reduce the problem to the study of deformations
of the base curve and of the fibers.
To this end, the first main ingredient is the fundamental result due
to Siu
and Beauville (cf. [B2], [C]) to the effect that the question of
whether a
surface admits a fibration over a curve of genus $\ge 2$ is of
completely
topological nature. In the same direction, fibrations obtained from
harmonic
maps were applied by Jost and Yau ([J-Y1], [J-Y2], [J-Y3]) to yield
related results which complement and strengthen the results obtained
by
algebro-geometric techniques. Combining these facts with the results
of Xiao on the structure of genus 2 fibrations ([X1]),
the study of the moduli spaces of the surfaces
under consideration readily reduces to
problems amenable to techniques from the theory of curves.

The case of smooth fibrations satisfying $q = g_{b} + 1$ is easier
to settle and in the first part of the paper we prove

{\bf Theorem 1} : {\it Given $g_{b} \ge 2$ and $K^{2}, \chi$
satisfying
$K^{2} = 8\chi = 8(g_{b} - 1)$. Then in the moduli space of surfaces
of
general type with these invariants $K^{2}, \chi$,
there exists a unique connected component of dimension
$3g_b - 1$
parametrizing smooth genus 2
fibrations with $q = g_{b} + 1$.}

In the rest of the paper, we deal with nonsmooth fibrations. In this
part
of the paper, the results of Xiao ([X1]) on universal genus 2
fibrations
and the work by Namba ([N]) on families of holomorphic maps from
compact Riemann surfaces into the projective line play crucial role.

We fix $1 \le K^{2}, \chi$ to satisfy $K^{2} < 8\chi$ and the
equality
$K^{2} = \lambda \chi(X) + (8 - \lambda)(g_{b} - 1)$
for some $g_{b} \ge 2$, where the slope $\lambda$ satisfies
$\lambda = 7 - \displaystyle \frac{6}{d}$,
$d \ge 2$. With this notation, we have the following partial result
which
complements results obtained by Seiler, because the surfaces we
consider do not admit genus 2 albanese fibrations.

{\bf Theorem 2} : \\
{\it (i) For $d \le 5$ and $g_b \ge 2$, depending on
$K^{2}, \chi$ there exists
a finite number of integers $n_i \ge 2$, i = 1,..,k
and for each $n_i$ a connected component of
dimension $\le 2n_i + 2g_b -4$ \\
parametrizing nonsmooth genus 2 fibrations with $q = g_b + 1$. \\
(ii) For $d \ge 6$ there exists a component parametrizing nonsmooth
genus
2 fibrations only if
$g_b \ge (d - 6)\displaystyle\frac{d^{2}}{24} \prod_{p \vert d}(1 -
\displaystyle\frac{1}{p^{2}}) + 1$.
For $d \ge 7$ and
$g_b = (d - 6)\displaystyle\frac{d^{2}}{24} \prod_{p \vert d}(1 -
\displaystyle\frac{1}{p^{2}}) + 1$,
we have a unique component of dimension 1.} \\

Throughout the paper we will consider only minimal surfaces over
${\bf C}$ and fibrations are always assumed to be connected. We adapt
the following standard notation.

%$\Delta = \{z \in {\bf C} : \vert z \vert < 1\}$.\\
%$\Delta^{*} = \Delta - \{0\}$.\\
$\chi(?), e(?), K_{?}, p_{g}(?), q(?)$ are the holomorphic Euler
characteristic,
the topological Euler characteristic,
the canonical class, the geometric genus and the irregularity of ?,
respectively. \\
$X$ is a smooth compact minimal surface admitting a fibration
$X \rightarrow S$ with base genus $g_{b} \ge 2$, fiber genus 2
and $q = g_{b} + 1$. \\
%$\Pi : {\cal X} \rightarrow \Delta$ is a degeneration of $X$, by
%%which we
%mean a proper flat holomorphic map from the smooth threefold ${\cal
%%X}$
%to the unit disk, such that $\Pi^{-1}(t_{0}) \cong X$ for
%some $t_{0} \in \Delta^{*}$ and $\Pi^{-1}(0)$ is the only singular
%%fiber. \\
$\pi_{1}(?)$ is the fundamental group of ?. \\
${\cal M}_{g}$ is the moduli space of curves of genus $g$. \\

We first observe that any deformation $X'$ of $X$ also admits a
fibration $X' \rightarrow S'$ with $q = g_b + 1$. In fact, since
$\pi_{1}(X') = \pi_{1}(X)$, by Beauville's criterion ([B2]) we have
a fibration $X' \rightarrow S'$ with $g(S') \ge g_{b}$. If $g(S') >
g_{b}$,
then as $g_b + 1 = q(X') \ge g(S')$ we have $g(S') = g_{b} + 1$ and
again by
Beauville's criterion we see that
$X$ is fibered by the albanese map too. This however is impossible,
because
then the first fibration on $X$ can not be connected. Therefore, to
understand the moduli of genus 2 fibrations with $q = g_{b} + 1$,
we must determine \\
(i) when the fibration $X' \rightarrow S'$ on a deformation $X'$ of
$X$
is of fiber genus 2, and \\
(ii) when two such genus 2 fibrations can be deformed to each other.

For smooth fibrations both of these questions can be easily answered.
In fact we will prove a general result for analytic fiber bundles of
base and fiber genera $\ge 2$. It is well known that such a fiber
bundle
with fiber $F$ becomes trivial over an unramified Galois base
extension
with group ${\cal G} = Image(\pi_{1}(S) \rightarrow Aut(F))$. In this
notation, we have

{\bf Proposition 1} : {\it Let $X \rightarrow S$ be a fiber bundle
of fiber genus $g \ge 2$ and base genus $g_b > (g + 1)/2$. Then each
deformation of $X$ admits an analytic fiber bundle
structure with the same trivializing group ${\cal G}$ as $X$.}

Proof :

We know by ([J-Y1], Lemma 7.1) that  fiber bundles
of genus $g$ with base of genus
$g_{b}$ deform to such bundles.
In what follows, we will prove that the trivializing group is fixed
in a family of fiber bundles parametrized by a connected base.

For a given ${\cal G}$, we take $g' = \vert {\cal G} \vert(g_{b} - 1)
+ 1$ and
we let $T_{g}, {\cal C}_{g}$ (resp. $T_{g'}, {\cal C}_{g'}$) be the
Teichm\"{u}ller
space of Riemann surfaces of genus $g$ (resp. $g'$) and the
universal curve on $T_{g}$
(resp. on $T_{g'}$). We consider the submanifolds
$T_{g}^{{\cal G}}$ of $T_{g}$
parametrizing curves $F$ with ${\cal G} \subset Aut(F)$
and $T_{g'}^{'}$ of $T_{g'}$ parametrizing curves of genus $g'$
admitting
${\cal G}$ as a fixed point free automorphism group.
We note that $dim(T_{g'}^{'}) = 3g_b - 3$, because ${\cal G}$ being
a quotient of $\pi_1(S)$, any curve of genus $g_b$ admits an etale
cover of genus $g'$.
We let ${\cal F} \rightarrow T_{g}^{{\cal G}} \times T_{g'}^{'}$
be the quotient by ${\cal G}$ of ${\cal C}_{g} \times {\cal C}_{g'}$
restricted to
$T_{g}^{{\cal G}} \times T_{g'}^{'}$.
For a point $p \in T_{g}^{{\cal G}} \times T_{g'}^{'}$ we consider
the
following diagram

\begin{tabular}{ccc}
$0 \rightarrow {\cal T}_{p}(T_{g}^{{\cal G}} \times T_{g'}^{'})$ &
$\rightarrow {\cal T}_{p}(T_{g} \times T_{g'})$ &
$\rightarrow {\cal T}_{p}(T_{g}^{{\cal G}} \times T_{g'}^{'})$ \\
&& \\
$\downarrow$ & $\downarrow$ & $\downarrow$ \\
&& \\
$0 \rightarrow H^{1}({\cal F}_{p}, \Theta)$ &
$\rightarrow H^{1}(F \times S, \Theta)$ &
$\rightarrow H^{1}({\cal F}_{p}, \Theta)$ \\
\end{tabular}

where the vertical arrows are the corresponding Kodaira-Spencer maps
and
the far right horizontal arrows are the trace maps induced by the
action
of ${\cal G}$. As the middle vertical arrow is an isomorphism, we
immediately
check that the Kodaira-Spencer map of the family
${\cal F} \rightarrow T_{g}^{{\cal G}} \times T_{g'}^{'}$
is bijective.
Hence this family is versal and
as $H^{0}({\cal F}_{p}, \Theta) = 0$ for all fibers ${\cal F}_{p}$
of the family, it is the universal deformation for any of its fibers.

Therefore, in a given deformation ${\cal X} \rightarrow {\cal S}$
of fiber bundles of genus 2, locally over
the base the trivializing group is fixed. Now, since
$g_b > (g + 1)/2$, we have
$K_{X}^{2} = 8(g_{b} - 1)(g - 1) > 4(g - 1)^{2}$ and therefore
none of the surfaces ${\cal X}_{s}$
admits two different genus g fibrations unless ${\cal X}_{s}$ is a
trivial product ([X1], Proposition 6.4). As trivial products deform
to
trivial products ([J-Y], Cor. 6.1), it follows that over a connected
base the trivializing group is fixed.$\Box$

{}From the proof of Proposition 1, it also follows that any two fiber
bundles
as in the statement of the proposition deform to each other.
Therefore, we see that in the moduli
space of surfaces with invariants $K^{2} = 8\chi = 8(g_b - 1)(g -
1)$,
analytic fiber bundles with fiber genus $g$ and base genus $g_b$, if
they exist, form a single connected component ${\cal M}$. The
existence for $g = 2$
case is proved in the following elementary lemma.

{\bf Lemma 2} : {\it Given $g_b \ge 2, K^{2} = 8(g_b - 1)$ and ${\cal
G}$ as
listed in ([X1], p. 30), there exists a genus 2 analytic fiber bundle
with
trivializing group ${\cal G}$.}

Proof :

Given ${\cal G}$, we take a curve $F$ of genus 2 with
${\cal G} \subset Aut(F)$.

We first observe that ${\cal G}$ is generated by two elements
$\{g_{1}, g_{2}\}$
if the hyperelliptic involution $\sigma$ on $F$ is not in ${\cal G}$
and by the four elements $\{g_{1}, g_{2}, \sigma g_{1}, \sigma
g_{2}\}$
if $\sigma \in {\cal G}$ where $\{g_{1}, g_{2}\}$ are coset
representatives
of $\langle \sigma \rangle$ in ${\cal G}$~ ([X1], p. 30). Therefore,
for a
curve $S$ of genus $\ge 2$ writing
$\pi_{1}(S) = \langle a_{j}, b_{j} : \prod [a_{j}, b_{j}] = e
\rangle$
we can define a surjective homomorphism

$\pi_{1}(S) \rightarrow {\cal G}$

by $a_{j} \mapsto g_{j}, b_{j} \mapsto g_{j}^{-1}$ (resp.
$a_{j} \mapsto g_{j}, b_{j} \mapsto \sigma g_{j}$) for $j = 1, 2$ if
$\sigma \notin {\cal G}$~ (resp. $\sigma \in {\cal G}$) and sending
all
other $a_{j}, b_{j}$ to identity. Via this homomorphism we construct
an
etale cover $S_{1} \rightarrow S$ with Galois group ${\cal G}$. The
quotient
surface $X = F \times S_{1}/{\cal G}$ admits a fibration $X
\rightarrow S$
of the type desired.$\Box$

For genus 2 fiber bundles with $q = g_b + 1$, ${\cal G} = {\bf Z}_2$
and the hyperelliptic involution $\sigma$ of the fiber is not in
${\cal G}$.
Taking ${\cal G}$ as described, the following lemma completes the
proof
of Theorem 1.

{\bf Lemma 3} : {\it Let ${\cal G} \neq \{e\}$ be one of the groups
listed
in ([X1], p. 30). Then
${\cal M} \cong {\cal M}_{2}^{{\cal G}} \times {\cal M}_{g_{b}}$,
where ${\cal M}_{2}^{{\cal G}}$
is the moduli space of curves of genus 2 admitting ${\cal G}$ as a
group of automorphisms.}

Proof :

Two fibers of the universal family ${\cal F}$ of the proof of
Proposition
1 are isomorphic if and only if they lie on $(x_{1}, y_{2}),
(x_{2}, y_{2}) \in T_{2}^{{\cal G}} \times T_{g'}^{'}$,
where $x_{1} \cong x_{2}$ (mod $\Gamma_{\cal G}$) and
$y_{1} \cong y_{2}$ (mod $\Gamma_{g_{b}}$), $\Gamma_{\cal G}$ being
the
modular group determined by $\cal G$. Hence
${\cal M} \cong T_{2}^{{\cal G}}/\Gamma_{\cal G} \times
T_{g'}^{'}/\Gamma_{g_{b}}
\cong {\cal M}_{2}^{{\cal G}} \times {\cal M}_{g_{b}}$. $\Box$ \\

For nonsmooth fibrations we first answer question (i) posed at the
beginning of the paper, for arbitrary fiber genus $g \ge 2$.

{\bf Lemma 4} : \\
{\it (i) A surface of general type $X$ admits, if any, a unique
nonsmooth fibration $X \rightarrow S$ with $g_b \ge g$
satisfying $q = g_{b} + 1$. \\
(ii) If $X$ admits such a fibration, then so does any
deformation $X'$ of $X$.}

Proof :

(i) By the criterion given in ([X1], Proposition 6.4),
it suffices to check that $K(X)^{2} > 4(g - 1)^{2}$.

For a fibration as in the statement of the lemma, we have

$K^{2} = \lambda \chi(X) + (8 - \lambda)(g_{b} - 1)$,

where $\lambda \ge 4$ ([X2], p. 459, Corollary 1). Therefore, as

$\chi(X) = (g_{b} - 1)(g - 1) -deg(R^{1}f_{*}({\cal O}_{X})) > (g_{b}
- 1)(g - 1)$

we get $K(X)^{2} > 4\chi(X) > 4(g - 1)^{2}$ as required.

(ii) Existence of the fibration on $X'$ is clear from the
discussion at the beginning of the paper. To see that fiber
genus is $g$, by the uniqueness of the fibration it suffices
to observe that the fiber genus is locally constant over the base of
a
deformation. But this is a consequence of ([\"{O}], Proposition 2)
which
in this context states that for a given deformation
${\cal X} \rightarrow {\cal Y}$
of $X$, locally around any $y \in {\cal Y}$ the deformation factors
over
a surface ${\cal C}$ to give a deformation of the fibers of
$X \rightarrow S$. $\Box$

As a consequence of Lemma 4, we see that in the moduli space of
surfaces
with invariants $K^{2}, \chi$, if there exist points corresponding to
surfaces
with nonsmooth genus 2 fibrations satisfying $q = g_b + 1$, then they
form connected
components. We will state our results on the existence and the
dimension
of such components with respect to the values of the invariant $d$ of
the
fibrations under consideration. First we recall the following basic
facts.

1) Given a genus 2 fibration $X \rightarrow S$ with invariant $d$,
the fibration is said to be of type $(E,d)$ if
the fixed part of the
associated jacobian fibration is the elliptic curve $E$.

2) For each pair $(E,d)$ where $E$ is an elliptic curve and $d \ge
2$,
there exists a genus 2 fibration $X(E, d) \rightarrow M(d)$ of type
$(E,d)$
over the modular curve $M(d)$ such that any fibration $X \rightarrow
S$
of this type is obtained as the minimal desingularization
of the pull back $f^{*}(X(E,d))$ of $X(E,d)$ via a surjective
morphism
$f : S \rightarrow M(d)$ ([X1], Corollaire on p.46).

Proof of Theorem 2 :

In the construction
of the universal fibration $X(E, d) \rightarrow M(d)$ given in
([X1], p.42), letting
$u_1(z) = (z_1, 0), u_6(z) = \displaystyle \frac{1}{d}(z_1, 1)$
vary, we obtain a family ${\cal X}(d) \rightarrow M(d) \times {\bf
A}^{1}$
where for $a \in {\bf A}^{1}$,
${\cal X}(d) \vert_{M(d) \times a} \cong X(E_a, d)$, $E_a$ being the
elliptic
curve with j-invariant $j(E_a) = a$.

Given a deformation ${\cal X} \rightarrow {\cal Y}$ of surfaces
admitting
fibrations as considered in this part of the paper, then around any
$y \in {\cal Y}$ the given deformation is induced from a deformation
of the
fibration ${\cal X}_y \rightarrow S_y$. Hence it follows that,
locally on
${\cal Y}$,
${\cal X} \rightarrow {\cal Y}$ is obtained via a
morphism $F$ into $M(d) \times {\bf A}^{1}$
where $F$ composed with the projection onto ${\bf A}^{1}$ maps $y$
to the j-invariant $j(E_y)$. Therefore, in the moduli
space of surfaces with fixed $K^{2}, \chi$, the dimension
of any component parametrizing surfaces with nonsmooth genus 2
fibrations
satsfying $q = g_b + 1$
is bounded by $m(d) + 1$ where $m(d)$ the dimension of the space of
holomorphic maps
from compact Riemann surfaces of genus $g_b$ into the modular curve
$M(d)$.

To prove part (i) of Theorem 2, we note that for
$d \le 5$, $M(d) \cong {\bf P}^{1}$ and we apply the results in ([N])
to get \\
a) in a deformation ${\cal X} \rightarrow {\cal Y}$ of $X$ over a
connected
base, the degree of the maps into $M(d) \cong {\bf P}^{1}$
inducing the fibration ${\cal X}_y \rightarrow S_y$ is independent
of $y$ ([N], Lemma 3.3.3), \\
b) for a fixed $n \ge 2$, the space of holomorphic maps of degree $n$
from compact Riemann surfaces of genus $g_b$ into
${\bf P}^{1}$ (modulo $Aut({\bf P}^{1})$) is of dimension $2n + 2g_b
- 5$
([N], Theorem 3.4.17).

On the otherhand, as a consequence of the restriction on $K^{2}(X),
\chi(X)$
if $X$ is obtained from $X(E,d)$ via a map of degree $n$ with
prescribed
ramification divisor,  it follows that for a given pair $K^{2}, \chi$
there exists, if any, only a finite number of integers $n_i \ge 2$
appearing as the degree of maps from Riemann surfaces of genus $g_b$
into $M(d)$ inducing surfaces with these given invariants. Combining
this
observation with (a) and (b) above, part (i) of Theorem 2 follows.

(ii) If there exists a nonsmooth genus 2 fibration of type $(E, d)$
over a
curve $S$ of genus $g_b$, then  we have a surjective map
$S \rightarrow M(d)$ and hence \\
$g_b \ge genus(M(d)) = (d - 6)\displaystyle\frac{d^{2}}{24} \prod_{p
\vert d}(1 - \displaystyle\frac{1}{p^{2}}) + 1$
for $d \ge 6$. \\
Furthermore, if we have $d \ge 7$ and $g_b = g(M(d))$, then since
$g(M(d)) \ge 2$, it follows that the map $S \rightarrow M(d)$ is an
isomorphism. Therefore, the moduli space of such fibrations
is ${\bf A}^{1}$, the modulus map being given by the j-invariant of
$E$ if $X$ is of type $(E, d)$. This completes the proof of Theorem
2. $\Box$ \\

{\bf References}

[B1] A. Beauville, L'in\'{e}galit\'{e} $p_{g} \ge 2q - 4$ pour les
surfaces de type
g\'{e}n\'{e}ral, Bull. Soc. Math. Fr. 110, 344 - 346 (1982).

[B2] A. Beauville, appendix to F. Catanese, Moduli and classification
of irregular K\"{a}hler manifolds (and algebraic varieties) with
albanese
general type fibration, Invent. Math. 104, 263 - 289 (1991).

[C] F. Catanese, Moduli and classification
of irregular K\"{a}hler manifolds (and algebraic varieties) with
albanese
general type fibration, Invent. Math. 104, 263 - 289 (1991).

[J-Y1] J. Jost, S. T. Yau, Harmonic mappings and K\"{a}hler
manifolds,
Math. Ann. 262, 145 - 166 (1983).

[J-Y2] J. Jost, S. T. Yau, Harmonic mappings and algebraic varieties
over function fields, Amer. J. Math. 115, 1197 - 1227 (1993).

[J-Y3] J. Jost, S. T. Yau, Harmonic maps and K\"{a}hler geometry,
in Prospects in Complex Geometry, LNM 1468, Sringer-Verlag (1991),
340 - 370 .

[\"{O}] H. \"{O}nsiper, A note on degenerations of fibered surfaces,
Indag. Math. 8(1), 115 - 118 (1997).

%[\"{O}-S] H. \"{O}nsiper, S. Sert\"{o}z, On the moduli spaces of
%fiber bundles of curves of genus 2, peprint.

[S1] W. K. Seiler, Moduli spaces for special surfaces of general
type,
in C. Bajaj (ed.), Algebraic geometry and its applications,
Springer-Verlag
(1994), 161 - 172.

[S2] W. K. Seiler, Moduli of surfaces of general type with a
fibration
by genus two curves, Math. Ann. 301, 771 - 812 (1995).

[X1] G. Xiao, Surfaces fibree en courbes de genre deux, LNM 1137,
Springer-Verlag (1985).

[X2] G. Xiao, Fibered algebraic surfaces with low slope, Math. Ann.
276,
449 - 466 (1987).

\vspace{0.5cm}
Department of Mathematics \\
Middle East Technical University \\
06531 Ankara, Turkey \\

e-mail : hursit@rorqual.cc.metu.edu.tr
\end{document}